\documentclass[12pt]{article}

\usepackage{amsfonts}  
\usepackage{graphicx}
\usepackage{pdfsync}

\title{ Overdetermined problems for fully non linear operators.} 

\author{I. Birindelli, F. Demengel}

\date{}

\catcode`@=11 \@addtoreset{equation}{section} \catcode`@=12

\newtheorem{theo}{Theorem}[section]
\newtheorem{prop}[theo]{Proposition}
\newtheorem{rema}[theo]{Remark}
\newtheorem{defi}[theo]{Definition}
\newtheorem{cor}[theo]{Corollary}
\newtheorem{lemme}[theo]{Lemma}

\def\R{{\rm I}\!{\rm  R}}
\def\grad{\nabla}

\def\bt{t^{\star}}

\begin{document}
\maketitle
\begin{abstract}
In this paper,  for $\alpha>-1 $, we consider the  overdetermined problem 
$ |\nabla u|^\alpha  {\cal  M}_{a,A} (D^2 u) = -f(u)$\ in a bounded  smooth domain  $\Omega$,  with Dirichlet condition 
 $u=0$  and Neumann condition $\partial_{\vec n} u = c \ {\rm on} \ \partial \Omega$ where $c$ is a constant,   $u$ is constant sign  and ${\cal M}_{a,A}$ is one of the Pucci's operator. We consider  different cases for  $f$, covering the case of the  principal eigenvalue for such operators. 
 In all the situations  considered  we  prove that, when $a$ is sufficiently close to $A$,  either $u= c=0= f(0)$, or $\Omega$ is  a ball, $u$ is radial, and $cu<0$ in $\Omega$. 
 \end{abstract}

\section{Introduction}
In this paper we prove that for a large class of nonlinearities $f(u)$, for ${\cal  M}_{a,A}$  one of the Pucci operators (i.e. either ${\cal  M}_{a,A}={\cal  M}_{a,A}^+$ or ${\cal  M}_{a,A}={\cal  M}_{a,A}^-$)  and $\alpha >-1$, if $\Omega$ is a bounded smooth domain,  such that there exists $u$ a 
viscosity,   constant sign  ${\cal C}^1$ solution of
 \begin{equation} \label{1} \left\{\begin{array}{lc}
 |\nabla u|^\alpha  {\cal  M}_{a,A} (D^2 u) +f(u)= 0&\ {\rm in} \ \Omega,\\
 u=0& \ {\rm on} \ \partial \Omega,\\
 \partial_{\vec n} u = c& {\rm on} \ \partial \Omega,
 \end{array} \right.
 \end{equation}
for some constant $c$, then 

$$\mbox{either}\quad c=f(0)=0\equiv u\quad \mbox{or}\quad \Omega\quad \mbox{is a ball and u is radial.} $$
Overdetermined boundary value problem is a very rich field, somehow started by the acclaimed paper by Serrin \cite{Se} where it is proved that, if $u$ is a solution of

$$ \left\{\begin{array}{lc}
 \Delta u = -1&\ {\rm in} \ \Omega\\
 u=0& \ {\rm on} \ \partial \Omega\\
 {\partial u\over \partial \vec n} = c& {\rm on} \ \partial \Omega,
 \end{array} \right.
$$
then $\Omega$ is a ball and $u$ is radial.
The proof relies on the method of moving planes. Let us remark that this method has already been extended to prove symmetry of solutions for fully nonlinear equations both by Gidas, Ni, Nirenberg \cite{GNNi} and by Da Lio, Sirakov \cite{DS}.

On the other hand the overdetermined problem has been greatly generalized to all kind of settings, and geometries and it would be far too long to enumerate all the interesting results achieved, let us remark that all these results concern divergence form operators. Instead, in order to motivate the results obtained here, we shall now  describe an interesting connection with principal eigenvalues. 

Precisely, let $\lambda(\Omega)$ be the functional that associates to a domain $\Omega$ the principal eigenvalue of the Dirichlet problem for the Laplace operator.
As it is well explained in \cite{EL} , a domain $\Omega$ is critical for the first eigenvalue functional under fixed volume variation if and only if the eigenfunction $\phi>0$ associated to $\lambda(\Omega)$ has constant Neumann boundary condition i.e. it is a solution of an overdetermined problem.
This is proved using the famous Hadamard equality (we refer to \cite{EL} and references therein).
In \cite{PSic}, Pacard and Sicbaldi have extended this result to Riemann manifolds.   

In recent years, the concept of principal eigenvalue has been extended to fully nonlinear operators, by means of the maximum principle (see \cite{BNV}).
The values
$$\lambda^+(\Omega)=\sup\{\lambda, \exists \phi>0 \ \mbox{in}\ \Omega,  |\nabla \phi|^\alpha  {\cal  M}_{a,A} (D^2 \phi)+\lambda \phi^{1+\alpha} \leq 0\quad\mbox{in}\ \Omega\}$$
$$\lambda^-(\Omega)=\sup\{\lambda, \exists \psi<0 \ \mbox{in}\ \Omega,  |\nabla \psi|^\alpha  {\cal  M}_{a,A} (D^2 \psi)+\lambda |\psi|^{\alpha}\psi \geq 0\quad\mbox{in}\ \Omega\}$$
are  generalized eigenvalues in the sense that there exists a non trivial solution to the Dirichlet problem 
$$|\nabla \phi|^\alpha  {\cal  M}_{a,A} (D^2 \phi)+\lambda^{\pm}(\Omega) |\phi|^{\alpha}\phi= 0\quad\mbox{in }\quad \Omega,\ \phi=0 \quad\mbox{on}\quad \partial\Omega.$$
 One of the open question, even for the Pucci operator is whether the Faber-Krahn inequality  holds or not in this context i.e. suppose that $\Omega$ is a domain of volume $V$ and suppose that $B$ is a ball with the same volume, is it true that
$$  \lambda^+(B)\leq \lambda^+(\Omega)?$$
A first step in this direction is to prove that the ball is critical for $\lambda^+(\Omega)$ under fixed volume variation. In view of what was described above for the Laplacian, the result obtained here i.e. that the only bounded domain for which the eigenfunction has constant boundary data is the ball, gives a good evidence that it may be the case that the ball is the only critical domain.

For unbounded domains the situation is slightly different, in \cite{Si}, B. Sirakov  considers the case of exterior domains and domains with several connected components and in this reference  he also proves that in order to have an overdetermined solution the domain has to be radial. Recently, in dimension 2, Helein, Hauswirth, and  Pacard in \cite{HHP} have constructed a domain for which there exists a harmonic function with zero Dirichlet data and  constant Neumann boundary,  which is neither radial nor an exterior domain. 
The construction of this domain is deeply related to the Laplace operator, but it would be interesting to know if a similar counterexample can be found for the Pucci operator. This will be the object of a future work.

We come now to a better description of the results contained in this note. It is well known that 
the last step in Serrin's proof is a sort of Hopf's lemma in "corners".
Indeed, if the domain contains a squared corner, and 
two ordered solutions touch each other at this corner, then, for any 
direction entering the domain, if the derivatives coincide then the 
second derivatives have to be separated. Interestingly, this result is  a
 consequence of the fact  that the eigenvalue of the Laplace Beltrami operator in a quarter sphere $S^{N-1}$ is  exactly $2N$, even though this is not obvious at all from Serrin's proof. In Proposition \ref{Hopfcorner} we extend Serrin's result to nonlinear setting considered here as long as $a$ is close to $A$.  Here
the difficulty is both that one needs to introduce a generalization of the Pucci's operator on the sphere and to estimate the eigenvalue on the quarter sphere; furthermore it is possible to prove that this eigenvalue is greater than $2N$. This is exactly the reason why we need to choose $a$ close to $A$.

The paper is organized in the following way, in the next section we state the results concerning the overdetermined problem, in the third section after recalling known results we prove a comparison principle which is new and interesting in itself, the last section is devoted to the proofs of the main result including the "Hopf lemma in corner" described above.

\section{The main result }
In the whole paper, for some $h\in(0,1)$, $\Omega$  is a  bounded  ${\cal C}^{2,h}$ domain of $\R^N$, $\alpha>-1 $, and $F$ is defined by
$$F(p, X) := |p|^\alpha {\cal M}_{a,A}(X)$$
 where either ${\cal M}_{a,A}={\cal M}_{a,A} ^+(X) =   A tr(X^+)-atr(X^-)$ or ${\cal M}_{a,A}={\cal M}_{a,A} ^-(X) =   a tr(X^+)-Atr(X^-)$. 
For  $f$   some  continuous  function we consider the overdetermined problem

  \begin{equation} \label{1} \left\{\begin{array}{lc}
|\nabla u|^\alpha  {\cal  M}_{a,A} (D^2 u) +f(u)=0&\ {\rm in} \ \Omega,\\
 u=0& \ {\rm on} \ \partial \Omega,\\
\partial_{\vec n} u= c& {\rm on} \ \partial \Omega
 \end{array} \right.
 \end{equation}
where $c$ is a constant and $\vec n$ denotes  the unit outer normal to $\partial \Omega$. 

 We shall consider the following three  cases:
   \begin{itemize}
\item[{\bf Case 1}] $f$ is nonincreasing and ${\cal C}^1$, $f(0)\geq 0$ .   
\item[{\bf Case 2}]  $f(u)=h(u)-g(u)$  with $h$ and $g$ odd,  continuous, non decreasing functions satisfying

$\forall s>1$, $\forall \tau>0$, $h(s\tau)\leq s^{1+\alpha}h(\tau)$, $g(s\tau)\geq s^{\beta}g(\tau)$ for some $\beta>1+\alpha$, and either $g>0$ on $\R^+$ or $g\equiv 0$. 

\item[{\bf Case 3}]  $\alpha = 0$,  $f$ is Lipschitz continuous. 
\end{itemize}

 \begin{theo}\label{th1} In these three cases, there exists a constant $\delta$ which  depends only on universal data and on $f$, such that for $|a-A|< \delta$, if there exists $u$  a constant sign ${\cal C}^1$ viscosity  solution of the overdetermined problem  (\ref{1}),  then
$$ \mbox{either}\quad c=f(0)=0\equiv u, \quad\mbox{or}\ \Omega \quad\mbox{is a ball},\quad u \ \mbox{is radial and}\quad u\ c<0.$$
 \end{theo}
 \begin{rema}
  In the case where $\alpha \leq 0$ the  ${\cal C}^1$ regularity of the solution is a consequence of the results in \cite{BD9,CC}. In the case where $\alpha >0$, except in the radial case,   in the one dimensional case or for  operators in divergence form,  this regularity is an  open  question. 
  \end{rema}
  
 \begin{rema}
  As an example in the case where 
  $f\equiv 1$ and for ${\cal  M}_{a,A}={\cal  M}_{a,A}^+$  one gets that the solution  is  given  by 
  $$\varphi (r)  = {\alpha+1\over \alpha+2} \left({1+\alpha \over a((N-1)(1+\alpha)+1)}\right)^{1\over 1+\alpha} \left(- r^{\alpha+2\over \alpha+1} +R_c^{\alpha+2\over \alpha+1}\right)$$ and $\Omega = B(0, R_c)$ where  $R_c$ and $c$ are linked by the relation 
    $$c =-  \left({1+\alpha \over a((N-1)(1+\alpha)+1)}\right)^{1\over 1+\alpha}  R_c^{1\over 1+\alpha}.$$  \end{rema}
As a consequence of Theorem \ref{th1}, in the case $f(u)=\lambda|u|^\alpha u$,   we get
\begin{cor}\label{th3} There exists a constant $\delta$ which  depends only on universal data, such that for $|a-A|< \delta$, 
the only bounded smooth domains for which an  eigenfunction with constant sign   satisfies
$$\partial_{\vec n} \psi=c\quad \mbox{on}\quad \partial\Omega$$
are balls.
   \end{cor}

   \begin{rema}
    The hypothesis that $a$ is close  to $A$  is only needed for the proof of Proposition  \ref{Hopfcorner} which is a generalization of the strict comparison in domains with corners in the case of the Laplacian, \cite{Se}.
    \end{rema}

\section{Preliminary results: comparison principles and regularity.}

   We begin by recalling   the definition of viscosity solution adapted to the  present context. 
  
  \begin{defi}\label{def1}
$v\in {\cal C} ({\Omega})\cap L^\infty(\Omega)$ is a viscosity super solution
of $ F(\nabla v, D^2 v)+ f(v) =0$ 
if, for all $x_o\in \Omega$, 

-either there exists an open ball $B(x_0,\delta)$, $\delta>0$  in $\Omega$
on which 
$v= cte= \kappa 
$ and 
$f(\kappa)\leq 0$, 

-or
 $\forall \varphi\in {\mathcal C}^2(\Omega)$, such that
$v-\varphi$ has a local minimum on $x_0$ and $\grad\varphi(x_0)\neq
0$, one has
\begin{equation}
F(\grad\varphi(x_0),
 D^2\varphi(x_0))+ f(v(x_0))\leq 0.
\label{eqdefi}\end{equation}
\end{defi}
Of course   a symmetric definition can be given  for the viscosity sub-solutions, and a viscosity solution is a function which is both a super-solution and a sub-solution.

We now recall some classical facts concerning the Pucci's  operators. 

\begin{prop}{\cite{CC}}\label{proppucci}
Suppose that $f$ is Lipschitz continuous and that  $u$ and $v$  are  respectively viscosity sub- and supersolutions of 
 $${\cal M}_{a,A} (D^2 w) + f(w)= 0\;\mbox{in}\ \Omega,$$ 
 and $u\leq v$ in $\Omega$. 
 
 Then either $u\equiv v$ or $u< v$ in $\Omega$ and 
 $\partial_{\vec n} (u-v) >0$ on $\partial \Omega$.
  \end{prop}
Furthermore a consequence of the famous Alexandrov-Bakelman-Pucci inequality allows to prove a maximum principle in "small domains":

\begin{prop}\label{small} Given $c(x)$ a bounded function in $\Omega$,  there exists $\delta$ depending on $|c|_\infty$ and on $a, A$, and the diameter of $\Omega$,  such that for any $\Omega_o\subset\Omega$ satisfying
$|\Omega_o|\leq \delta$:
$$\left\{\begin{array}{lc}
{\cal M}_{a,A}(D^2 w) + c(x)w\geq 0&\mbox{in}\ \Omega_o,\\
w\leq 0 &\mbox{on}\ \partial\Omega_o
\end{array}\right. \Rightarrow w\leq 0\; \mbox{in}\ \Omega_o.$$ 
\end{prop}
The proof is well known (see \cite{BNV}) but we recall it for completeness sake.
Observe that $w$ satisfies
$${\cal M}_{a,A} (D^2 w) - c^-(x)w\geq -c^+ w \ \mbox{in}\ \Omega_o.$$
Hence  the Alexandroff Backelman Pucci's theorem implies that $w$  satisfies (\cite{CC})
$$\sup_{\Omega_o} w\leq C\|c^+ w\|_{L^N(\Omega_o)}$$
where $C$ is a constant that depends on $a$, $A$ and the diameter of $\Omega$.
Hence for $|\Omega_o|$ sufficiently small, $\sup_{\Omega_o}w\leq 0$.

We shall also need the following regularity result in the case $\alpha = 0$, \cite{Vi}, \cite{E,CC}.

 \begin{prop}\label{rglr1}   Let $f$ be some  bounded and H\"older function on $\Omega$.
   Then for all $A>0$ there exist $\kappa=\kappa(A, f, \Omega)$ and $C=C(A, f, \Omega)$,  there exists $\epsilon >0$ such that for all $t\in ]1-\epsilon,  1], $  and any $u$ viscosity solution of 
  $$\left\{ \begin{array}{lc}
  {\cal M}_{tA,A} (D^2 u) = f & {\rm in} \ \Omega\\
  u= 0 &{\rm on}\  \partial \Omega, 
  \end{array}\right.$$
   satisfies 
   $$||u||_{{\cal C}^{2, \kappa}(\overline{\Omega})} \leq C. $$
\end{prop}
This will also be used in the case $\alpha\neq 0$. 

Comparison principles play a key role when one deals with viscosity solutions.
We both recall known one (Theorem \ref{comp})
 and prove a new one (Theorem \ref{comp3}).

  \begin{theo}{\cite{BD1}}\label{comp} Suppose that
$\phi$ and $\sigma$ are respectively,  sub- and super-solutions of 

$$F( \grad \phi,D^2\phi)-\beta (\phi)\leq f_1 \ {\mbox in} \ \Omega, $$
$$F(  \grad \sigma,D^2\sigma) -\beta (\sigma)\geq f_2\ {\mbox in} \ \Omega,$$ 
with $f_1$,  $f_2$ and $\beta $  continuous functions  on $\R^+$ such that 
 
 -either $\beta$ is increasing on $\R^+$ and $f_1\leq f_2$, 
 
 -or $\beta$ is
nondecreasing and $f_1< f_2$. 

\noindent If $\sigma\leq\phi$ on $\partial
\Omega$ then $\sigma \leq \phi$ in $\Omega$. 
\end{theo}
 For the proof of Theorem \ref{th1}  we shall need the  following refined  comparison principle, where we have denoted in a classical way and for simplicity $F[v] = F(\nabla v, D^2 v)$ : 
           
\begin{theo}\label{comp3} Assume that  $u$ and $v$ are constant sign , $|v|>0$  on $\overline{\Omega} $, and   are   viscosity solutions of 
$$F[v] + h(v)-g(v) \leq 0\; \mbox{in} \ \Omega$$
             and 
 $$F[u] + h(u)-g(u) \geq 0\;  \mbox{in} \ \Omega, $$
where  $h$ and $g$ are continuous, odd and  non decreasing functions such that for some $\beta >1+\alpha $,  for all $s>1$ and for all $\tau>0$

\begin{itemize}
\item  $h(s \tau) \leq s^{1+\alpha} h(\tau)$,  
   
   \item  $g(s\tau) \geq s^\beta g(\tau)>0$.
           
\end{itemize}
Then the comparison principle holds i.e. if $u\leq v$ on $\partial \Omega$ then
$u\leq v$ in $\Omega$. 

If $g\equiv 0$ and $h$ is increasing  then the same conclusion holds.

\end{theo}
The proof is postponed to the end of the section.
   \begin{rema} In these Theorems,  $\Omega $  needs not   be  regular, bounded is sufficient. 
\end{rema}

 We shall also need the following  strong  comparison principle  : 
  
 \begin{prop}\cite{BD9}\label{chopf}
Let $f$ be  ${\cal C}^1$ and let  $u$ and $v$ be respectively nonnegative 
${\mathcal C}^{1}(\overline{\Omega})$ viscosity  solutions of

$$ 
 F(\nabla u,D^2 u)+ f(u)\leq 0\  \mbox{in}\ \quad\Omega,
 $$
$$
 F(\nabla v,D^2 v)+  f(v)\geq 0\  \mbox{in}\quad\Omega,
 $$
 with $u\geq v$ in $\Omega$. 
Suppose that there exists $\bar x\in \overline{\Omega}$  such that $u(\bar x) = v(\bar x)$ and either $\nabla v(\bar x)\neq 0$ or $\nabla u(\bar x)\neq 0$,
then  there exists $R$ such that $u\equiv v$ on $B(\bar x, R)\cap \Omega$. 

Furthermore  if $v>0$ in $\Omega$,  $v=0$ on $\partial \Omega$,   such that  there exists $\bar x\in\partial\Omega$  such that  $u(\bar x)=0$,  and $\partial_{\vec n} u(\bar x)= \partial_{\vec n} v(\bar x)$,   
 then there exists $\epsilon>0$ such that 
$$u\equiv v\ \mbox{in}\ \Omega\setminus\overline\Omega_\epsilon$$
where $\Omega_\epsilon$ is the set of points of $\Omega$ whose distance to  the  connected component of the boundary  which contains $\bar x$ is greater than $\epsilon$.
\end{prop}
This proposition holds for a more general class  of   operators  than the one considered here. It will  be used in the proof of Theorem \ref{th1}.

{\em Proof  of Theorem  \ref{comp3}.} We can assume without loss of generality that $u$ and $v$ are positive. 
 
 We suppose by contradiction that somewhere $u> v$.  Let $\gamma^\prime = \sup_{\Omega}  {u\over v}$,   let 
 $\kappa =\displaystyle \inf_{x\in \overline{\Omega}} |g(v(x))|((\gamma^\prime) ^{\beta }-(\gamma^\prime )^{1+\alpha})$
  and let $\gamma\in ]1, \gamma^\prime[$ sufficiently close to $\gamma$ in order that 
  $$\sup_{x\in \overline{\Omega}} |h(\gamma v)-h(\gamma^\prime v)|\leq {\kappa\over 4}
  $$
and  $ \displaystyle\inf_{x\in \overline{\Omega}} |g(v(x))|(\gamma^{\beta }-\gamma ^{1+\alpha})\geq {3\kappa\over 4}$. 
            Let us note that  $u-\gamma v$ achieves its positive maximum  inside $\Omega$.     
                
 Let us define   $\psi_j(x,y) = u(x)-\gamma  v(y)-{j\over q} |x-y|^q$ where 
 $q> \sup ({\alpha+2\over \alpha+1}, 2)$. 
It is classical that $\psi_j$ achieves its maximum on some pair $(x_j, y_j)$ which is in $\Omega^2$ and that $(x_j, y_j)\rightarrow (\bar x, \bar x)$ where 
$$u(\bar x)-\gamma v(\bar x)=\sup_{x\in \overline{\Omega}}( u( x)-\gamma v( x))>0.$$ 
Moreover $j|x_j-y_j|^q\rightarrow 0$. Then using Ishii's lemma, see  \cite {I,BD1} , there exist  $X_j$, $Y_j$ in $S$ with   $(j|x_j-y_j|^{q-2} (x_j-y_j), X_j)\in J^{2,+} u(x_j)$, $(j|x_j-y_j|^{q-2} (x_j-y_j), -Y_j)\in J^{2,-}  v(y_j)$ with 
$$X_j+ \gamma Y_j\leq 0.$$
In order to use the equations, from the definition of viscosity solutions we need to prove that   $x_j\neq y_j$ , this  will be checked  later.  One has, using the fact that $u$ and $v$ are sub and super solutions
 \begin{eqnarray}\label{eqetoile}
                -h(u(x_j) ) + g(u(x_j))&\leq& F(j|x_j-y_j|^{q-2} (x_j-y_j), X_j) \nonumber \\
                 &\leq&  \gamma^{1+\alpha} F(j|x_j-y_j|^{q-2} (x_j-y_j), -Y_j)\nonumber\\
                 &\leq &  \gamma^{1+\alpha}\left(-h(v(y_j))+g(v(y_j))\right).\nonumber\\
\end{eqnarray}
Passing to the limit  and  using the properties of $h$ and $g$ one obtains 
\begin{eqnarray*}
  -h(\gamma^\prime v(\bar x) )+ g(\gamma  v(\bar x) )&\leq &  -h(u(\bar x))+ g(u(\bar x) )\\
  &\leq& \gamma^{1+\alpha} (-h(v(\bar x))+ g(v(\bar x)))\\
   &\leq & -h(\gamma v(\bar x))+ \gamma^\beta  g(v(\bar x))-{3\kappa\over 4}\\
   &\leq & -h(\gamma ^\prime v(\bar x) )+ g(\gamma v(\bar x))-{\kappa\over 2}
   \end{eqnarray*}
which is a contradiction.

   \noindent We now suppose that $g\equiv 0$ and $h$ is increasing.  We begin to prove the result when there exists $\delta>0$ such that
\begin{equation}\label{eq30}F[v] + h(v)\leq -\delta
      \end{equation}
Since $v>0$ on $\overline{\Omega}$,  we define $\gamma^\prime $ as before,
we want to prove that $\gamma^\prime\leq 1$, then we suppose by 
contradiction that $\gamma^\prime  >1$.  Let $\gamma \in ]1, \gamma^\prime[$
be small enough in order that by the continuity of $h$  and the boundedness of 
$v$ one has

$$\sup_{x\in \overline{\Omega}}|h(\gamma^\prime v(x))-h(\gamma v(x))|\leq {\delta\over 4}.  $$ 
By passing to the limit  in (\ref{eqetoile}) with $g\equiv 0$,  and using the properties of  $h$, we obtain 

$$-h(\gamma^\prime  v(\bar x))\leq -h(u(\bar x) )\leq -\gamma^{1+\alpha}h(v(\bar x))-\delta \leq -h(\gamma^\prime v(\bar x))-\frac{\delta}{2},$$
a contradiction.

Suppose (\ref{eq30}) does not hold, and recall that $v>0$ on 
$\overline{\Omega}$. For any arbitrary positive $\epsilon$ let $w_\epsilon = v(1+\epsilon)-{\min_{\Omega} v\over 2} \epsilon$. Then   $u< w_\epsilon$ on $\partial \Omega$ and since $h$ is  now  supposed to be increasing,  there exists $\delta_\epsilon>0$ such that $h(w_\epsilon)\leq (1+\epsilon)^{1+\alpha} h(v)-\delta_\epsilon  $ hence 

$$F[w_\epsilon]+ h(w_\epsilon) \leq (1+\epsilon)^{1+\alpha} (F[v]+ h(v))-\delta_\epsilon \leq -\delta_\epsilon$$
and then, from the previous result, $u\leq w_\epsilon$ in $\Omega$ and, letting $\epsilon$ go to zero, $u\leq v$ in $\Omega$.

There remains to prove that  $x_j\neq y_j$ definitively. If $x_j= y_j$, one would have 
$$v(x) \geq v(x_j)-{j\over q} |x-x_j|^q\quad \mbox{and}\quad u(x) \leq u(x_j)+ {j\over q} |x-x_j|^q.$$
If the infimum 
$$\inf_{x\in B_r(x_j)} \{ v(x)+ {j\over q} |x-x_j|^q\}$$ is not strict then one can replace $x_j$ by some point $y_j$ close to it and then we are done.  The same is true if we assume that the supremum 
$$\sup_{x\in B_r(x_j)}\{  u(x)-{j\over q} |x-x_j|^q\}$$ 
is not strict. So we assume that  both extrema are strict. In this case, proceeding as in \cite{BD1} one can prove, using the equation and the definition of viscosity solution, that
                 
$$h(v(x_j))-g(v(x_j)) \leq 0\quad \mbox{and}\quad h(u(x_j))-g(v(x_j))\geq 0.$$ 
Passing to the limit the inequality becomes
$$h(v(\bar x))-g(v(\bar x))\leq 0\quad \mbox{and}\quad h(u(\bar x)) -g( u(\bar x)) \geq 0.$$
                   Using $u(\bar x) > v(\bar x)$ one derives that 
$$h(u(\bar x)) \geq g(u(\bar x)) \geq\left( {u(\bar x)\over v(\bar x)}\right)^\beta g(v(\bar x)) \geq  \left({u(\bar x)\over v(\bar x)}\right)^\beta  h(v(\bar x)).$$
Let us note that we have $h(u(\bar x)) \geq 0$ by the previous inequalities and then also 
 $$h(v(\bar x)) \geq \left({v(\bar x)\over u(\bar x)}\right) ^{1+\alpha} h(u(\bar x)) >0.$$
 Finally this gives 
$$\left(\left({u(\bar x)\over v(\bar x)}\right)^{1+\alpha}-\left({u(\bar x)\over v(\bar x)}\right)^\beta \right) h(v(\bar x)) >0$$
which is a contradiction, since ${u(\bar x)\over v(\bar x)}>1$ and $\beta > 1+\alpha$.   
                    
In the case where $g\equiv 0$ the result holds by the increasing behavior of $h$.  This ends the proof of Theorem \ref{comp3}.

  \bigskip   
\noindent We end this section with an important remark concerning regularity of solutions close to the boundary :

\begin{rema}\label{rglr} Observe that, as a consequence of Proposition \ref{rglr1}, using Hopf lemma, we know that  for any $u$, ${\cal C}^1$, constant sign solution of
  $$ \left\{ \begin{array}{lc}
 |\nabla u|^\alpha \ {\cal M}_{a,A}(D^2 u) +  f (u)=0& \mbox{in} \ \Omega, \\
u = 0 & \mbox{ on} \ \partial \Omega,
\end{array}\right.$$
 then there exists $\gamma\in (0,1)$ and a neighborhood of $\partial\Omega$ such that $u\in {\cal C}^{2,\gamma}$ in that neighborhood.
 \end{rema}
To prove this regularity in the case $\alpha <0$, this hypothesis that $u$ is ${\cal C}^1$ is not needed, furthermore the result is true everywhere; the proof can be found in \cite{BD9}. When $\alpha  >0$ one can use the same arguments as in \cite{BD11},   Theorem 2.8.

\section{Proofs of the main results}

As in Serrin's  original paper  \cite{Se} we use the moving planes method.

  We shall need the two following results :

   \begin{prop}\label{Hopfcorner}
   Suppose that $f$ is ${\cal C}^1$ on $\R^+$. 
    Suppose that  $\Omega^\star$ is  some bounded ${\cal C}^{2, h}$ domain, and suppose that $H_0$ is an hyperplane such that there exists $P\in H_0\cap \partial \Omega^\star$, with $\vec n_{\Omega^\star} (P)\in H_0$. Let $\Omega$ be  the intersection of  $\Omega^\star$ with one of the half spaces bounded by $H_0$. 
    
Suppose that $u$ and $v$ are ${\cal C}^2$ solutions of
$$\left\{ \begin{array}{lc}
    |\nabla v|^\alpha {\cal M}_{a,A} ( D^2 v)+ f(v)\leq  |\nabla u|^\alpha {\cal M}_{a,A} ( D^2 u)+ f(u) \  \mbox{ in}
     \ \Omega,\\
      u<  v\  \mbox{ in a  neighborhood of } \ P\ \mbox{ in } \ \Omega,   &  \  \\
      u(P) = v(P) \  \mbox{ and either} \ |\nabla u(P)|\neq 0 \ \mbox{ or} \  |\nabla v(P) | \neq 0.& 
      \end{array}\right. 
 $$
 For any $\vec \nu \in \R^N$  a direction pointing inside $\Omega$  i.e.  such that  $\vec \nu\cdot \vec n(P) <0$, and also such that $\vec \nu\cdot e_1 >0$, then,       
$$\mbox{either}\quad \partial_{\vec \nu} v (P) > \partial_{\vec \nu} u(P)\quad \mbox{or} \quad\partial^2_{\vec \nu} v(P) > \partial^2_{\vec \nu}  (P) .$$
 \end{prop}
 
\begin{lemme}\label{lem1}
For any $u$ solution of (\ref{1}), if $\partial\Omega$ is the zero level set of a
function $\psi$ then for any $P\in \partial \Omega$,   
$D^2 u(P)$ depends only on $\psi$, $\nabla \psi$ and $D^2\psi$ on $P$.         
 \end{lemme}
 
  We postpone the proofs of these two results and prove 
 Theorem  \ref{th1}. For convenience of the reader we recall the three cases we are going to treat:
  \begin{itemize}
\item[{\bf Case 1}] $f$ is nonincreasing  and ${\cal C}^1$, $f(0)\geq 0$.   
\item[{\bf Case 2}]  $f(u)=h(u)-g(u)$  with $h$ and $g$ odd,  continuous, non decreasing functions satisfying

$\forall s>1$, $\forall \tau>0$, $h(s\tau)\leq s^{1+\alpha}h(\tau)$, $g(s\tau)\geq s^{\beta}g(\tau)$ for some $\beta>1+\alpha$, and either $g>0$ on $\R^+$ or $g\equiv 0$. 

\item[{\bf Case 3}]  $\alpha = 0$,  $f$ is Lipschitz continuous. 
\end{itemize}

 \bigskip
  { \bf Proof of Theorem \ref{th1}} 
We start by remarking that by Hopf's principle either $u\equiv 0$ and then $c= f(0) = 0$, or $|u|>0$ in $\Omega$. Without loss of generality we shall suppose that $u>0$ and then $c<0$.

In order to start the moving plane procedure, we choose a direction, say $e_1$, and
for $t\in \R$, we denote by   $H_t $   the hyperplane $\{x_1=t\}$ and   the sets 
  $\Omega_t^- = \Omega \cap \{ x_1<t\}$, and $\Omega^+_t =\{ x, x_1>t, (2t-x_1, x^\prime )\in \Omega_t^-\}$.  
  
We  define $u_t (x) = u(2t-x_1, x^\prime)$. It is easy to see that for any $\phi\in{\cal C}^2$, 
the eigenvalues of the Hessian of $\phi$ and $\phi_t$ are the same, hence, using the definition of viscosity solution and the definition of Pucci's operator, we get
that $u$ and $u_t$ satisfy the same equation in $\Omega_t^+$.

   It is clear that for $t<0$ large, $\Omega_t^- = \emptyset$.   
Let $t_1 = \sup\{ t,   \Omega_t^-  = \emptyset\}$ and $\bt = \sup \{ \tilde t,\ \forall t<\tilde t\ , \  \Omega_t^+ \subset \Omega\}$ then $\bt$ is such that one of the two following events  occurs: 
  
  - {\bf event 1}  : $ H_{\bt}$ contains the normal to the boundary of $\Omega$ at some point $P$, or  
  
   -  {\bf event 2} : $\Omega_{\bt}^+$ becomes internally tangent to the boundary of $\Omega$ at some point $P$ not on $H_{\bt}$. 

Recall that for any $t\in (t_1,\bt)$,  $u=u_{t}$ on  $H_{t} \cap \Omega_{t}^+$,  and $u\geq u_{t}$ 
on
    $\partial  \Omega_{t}^+\cap \Omega$.
    
In all three cases we need to prove the following two steps:

{\bf Step 1} $u_t \leq u$ in $\Omega_t^+$ for any $t   \in (t_1,\bt]$.

{\bf Step 2}  $\Omega$ is symmetric with respect to $H_{\bt}$ i.e. $\Omega=\Omega_{\bt}^-\cup \Omega_{\bt}^+\cup H_{\bt}$.

This ends the proof because  since the direction $e_1$ was chosen arbitrarily,  this implies that $\Omega$ is symmetric with respect to any direction and is therefore a ball.

\bigskip
{\bf Proof of step 2.}  First suppose that "event 2" occurs i.e. there exists $P\in \partial \Omega_{\bt}^+ \cap  \partial \Omega$. 
Since the unit exterior normal to $\partial\Omega$ in $P$ is the same than the one of $\partial\Omega-{\bt}^+$ and by obvious symmetries,
$\partial_{\vec n}u_{\bt} (P)= \partial_{\vec n}u (P)=c$.  Using  Proposition \ref{chopf}, one gets that $u =  0$ on all the connected component    of $\partial \Omega_{\bt}^+\cap \overline{\Omega}$ which contains $P$,  this implies that $\partial \Omega \cap \partial \Omega_{\bt}^+= \partial \Omega_{\bt}^+ \setminus H_{\bt}$. Then $\Omega$ is symmetric with respect to $H_{\bt}$, hence symmetric with respect to $x_1$.

We  now consider "event 1", i.e. we suppose that  there exists some point $P\in H_{\bt}\cap \partial \Omega$, with $\vec n_\Omega (P)\in H_{\bt}$. We begin to prove that $u = u_{\bt}$ in a neighborhood of $P$ in $\Omega_t$. 

 Since $\nabla u \neq 0$  around $P$, using Proposition \ref{chopf} either $u\equiv u_{\bt}$ or  $u> u_{\bt }$ 
in a neighborhood of $P$. 
 
Suppose by contradiction that $u > u_{\bt}$    inside $\Omega_{\bt}\cap B(P, R)$, then, by Proposition \ref{Hopfcorner}, if $\vec\nu$ is such that $\vec\nu\cdot\vec n <0$, and $\vec\nu\cdot \vec  e_1>0$,  either $\partial_{\vec \nu} u (P) > 
\partial _{\vec \nu} u_{\bt} (P)$ or $\partial^2_{\vec \nu} u (P) > \partial^2 _{\vec \nu} u_{\bt} (P)$. 

The first inequality is impossible since  on $\partial \Omega$, $\partial_{\vec \nu} u (P) = c (\vec \nu\cdot \vec n)= \partial _{\vec \nu} u_{\bt} (P) $.  The second inequality is also impossible because Lemma \ref{lem1} implies 
that $\partial_{\vec \nu}^2 u(P) = \partial_{\vec \nu}^2 u_{\bt} (P)$.

Observe that in case 2  one  applies  Proposition  \ref{Hopfcorner}  in the following manner  
$$F[u] -g(u)\leq -h(u) \leq -h(v) \leq F[v]-g(v).$$

We have obtained that  $u= u_{\bt}$ in a neighborhood of $P$. This implies in particular that $u=0$ on $\partial \Omega_{\bt}\cap B(P,R)$  hence $\partial \Omega_{\bt} \cap B(P, R) \subset  \partial \Omega$. Using Proposition \ref{chopf} we 
get that $u=0$ in $\partial \Omega_{\bt}\setminus H_{\bt}$.
This of course implies that $\Omega$ is symmetric with respect to $H_{\bt}$.

\bigskip
{\bf Proof of Step 1 in Case 1} is just an application of Theorem  \ref{comp} in $\Omega_t^+$.

\bigskip
{\bf Proof of Step 1 in Case 2.}
For $t < \bt$  there are no points in $\partial \Omega \cap H_t$ with $\vec n_\Omega \in H_t$. Then, for $\vec \nu=-{\vec n}_{\Omega_t} (P)$,
$$\partial_{\vec \nu} u (\bar x)>0\ \mbox{ and }\ \partial_{\vec \nu} u_t (\bar x)<0$$
 As a consequence  there exists $\epsilon >0$ such that on $B(\bar x, \epsilon) \cap \Omega_t^+$,  $u_t\leq  u$.  Let $B_{\epsilon}=\cup_{\bar x\in \partial \Omega \cap H_t} B(\bar x, \epsilon)$. 
 
 Since $u>0$ on $\partial \Omega_t^+\setminus B_\epsilon$ there exists a neighborhood $V_\epsilon$ 
 of $\partial \Omega^+_t  \setminus B_\epsilon $ such that $u> u_t$ on $V_\epsilon$. Let $\Omega_o = \Omega_t^+\setminus (B_\epsilon\cup V_\epsilon) $,  then $u_t>0$ in $\bar\Omega_o$ and $u> u_t$ in $B_\epsilon \cup V_\epsilon$.  

 We are in the hypothesis of Theorem \ref{comp3} hence $u_t\leq u$ in $\Omega_o$ and hence in $\Omega_t^+$. 
By continuity, the inequality holds also for $t=\bt$.

\bigskip
{\bf Proof of Step 1 in Case 3.}

Let us recall that we are in the case $\alpha = 0$, and $f$ is only supposed to be Lipschitz continuous.  Here the  key argument will  not be the first comparison principle in Theorem \ref{comp} but  the maximum principle in small domains.

  We start by proving that, for $t$ sufficiently close to $t_1$, $u_t\leq u$ in $\Omega_t^+$. 
      
\noindent Without loss of generality one can assume that  $t_1 = 0$. We need to prove that for some $h>0$ and for $t\in [0, h[$, $u_t \leq u$ in $\Omega_t^+$. 
      
       \noindent Let $Q\in \partial \Omega \cap H_0$. Then $\vec n_\Omega (Q) = -e_1$,   Neumann condition implies that  
       $\partial_{x_1} u (Q)  =- c$, hence since $u$ is ${\cal C}^1$,  there exists $r>0$ such that on $B(Q, r) \cap \Omega$,  $\partial_{x_1} u(x) \geq {-c\over 2}$. Hence for $t$ small enough $u$ is strictly increasing in $\Omega\cap\{x_1<2t\}$ and then, for $t<x_1<2t$, 
   $u(2t-x_1, x^\prime) < u(x_1, x^\prime)$. 
        
\noindent We now define 
 $$\bar t = \sup \{ t\leq t^\star, \forall t^\prime < t,\ u_{t^\prime} \leq u\ {\rm in} \ \Omega_{t^\prime}^+\}.$$
We want to prove that $\bar t = t^\star$. 
            
Suppose by contradiction that $\bar t< t^\star$ then $\partial \Omega_{\bar t}^+\setminus H_{\bar t} \subset \Omega$ with $n_\Omega (P)\cdot \vec e_1 <0$ on $P\in \partial \Omega \cap H_{\bar t}$. These  two  conditions imply that  $\Omega_{\bar t+ h} \subset \Omega$ for $h$ small enough.  
Observe that $u_{\bar t}< u$ in $\Omega_{\bar t}^+$. Indeed,  since $f$ is Lipschitz continuous, one can use the strong maximum principle Proposition \ref{proppucci}  for the difference $u_{\bar t}-u$ and obtain both that $u_{\bar t}< u$ inside $\Omega_{\bar t}^+$
 and $\partial_{x_1} (u-u_{\bar t}) >0$ on $\partial \Omega_{\bar t}\cap  H_{\bar t}$. 
   
   {\bf Claim} For $h>0$ small enough $u_{\bar t+ h} \leq u$ in $\Omega_{\bar t+ h}^+ $. 
   
  \noindent  This claim  will contradict the definition of $\bar t$. 
    
\noindent To prove  the claim, let $K$ be a  compact subset of $\Omega^+_{\bar t}$ such that
$$|\Omega^+_{\bar t}\setminus K|\leq 2\delta,$$
where $\delta>0$ is the constant in Proposition \ref{small} with respect to $\Omega$ and $|\gamma(x)|_\infty=L_f$ the Lipschitz constant of $f$ .Clearly in $K$, $u_{\bar t}<u$ and, by continuity, for any $h$ sufficiently small, we still have $u_{\bar t+h}<u$ in $K$.

\noindent Take  $h$ sufficiently small  in order  that  $K\subset \subset \Omega_{\bar t+ h}$ and 
$$|\Omega^+_{\bar t+h}\setminus K|\leq \delta.$$
Since $u$ and $u_{\bar t+h}$ satisfy the same  equation in $\Omega_{\bar t+ h}$, 
$w=u_{\bar t+h}-u$  satisfies 
$${\cal M}_{a,A}^+ w+L_f w\geq 0\ \mbox{ in }\ \Omega^+_{\bar t+h}\setminus K$$
and  $w\leq 0$ in $\partial \left(\Omega^+_{\bar t+h}\setminus K\right)$.

Applying Proposition \ref{small}, we obtain that $w\leq0$ in  $\Omega^+_{\bar t+h}\setminus K$.
Finally $u_{\bar t+h}\leq u$ in $\Omega_{\bar t+h}$, for any $h>0$ sufficiently small.

 \noindent    We have obtained that $\bar t = t^\star$.
 This ends the proof of  step 1 and hence of Theorem \ref{th1}.

\bigskip
We now prove the two technical results used in the proof.

The proof of Proposition \ref{Hopfcorner} relies on a more general lemma about barriers by below on  some angular sector in the sphere 

\begin{lemme}\label{lemw}
 For any $S$ an open connected subset of the quarter  sphere $S^+ $,  for any $\epsilon >0$, there exists $\gamma_{a,A,\epsilon}^S>0$,  and $\psi : S\rightarrow \R$ such that for any $\gamma
  \geq \gamma_{a,A,\epsilon }^S$,  the function $w = r^{\gamma} \psi_S (\sigma)$
satisfies
 $$\left\{\begin{array}{lc}
{\cal M}^-_{a,A} (D^2 w) \geq  \epsilon r^{-2}  w & {\rm in}\ \{x\in\R^N;\ \frac{x}{|x|}\in S\},\\
 \psi >0 \ &\ {\rm in} \ S,\\
 \psi = 0 \ & {\rm on} \ \partial S
 \end{array}\right.$$ 
  and  $|\nabla w|\leq C_{a,A, S, \gamma}  r^{\gamma-1}$, for some constant  $C_{a,A, S, \gamma} $. 
 
 Furthermore if $S_\delta\subset \subset S^+$, with $|S_\delta -S^+|\leq \delta$, 
 $\lim_{(\epsilon, \delta) \rightarrow (0,0), \  a\rightarrow A} \gamma_{a,A, \epsilon}^{S_\delta} = 2$.
   \end{lemme}

We postpone the proof of Lemma \ref{lemw} and prove Proposition \ref{Hopfcorner}. 

\bigskip
\noindent{\it Proof  of Proposition \ref{Hopfcorner}.} Without loss of generality we shall suppose that $H_o=\{x_1=0\}$.
Let us note first that since $v>u$ on a neighborhood of $P$,  $\partial_{\vec \nu} v(P)\geq \partial_{\vec \nu} u(P)$, so we assume that $\partial_{\vec \nu} v(P) = \partial_{\vec \nu} u(P)$, and we want to prove that 
$\partial^2_{\vec \nu}  v(P) > \partial^2_{\vec \nu} u (P) .$

Since the boundary of $\partial \Omega^\star$ is ${\cal C}^{2,h}$  and $u$ is ${\cal C}^1$, one can assume that $0=P$ and  $R$  is such that 
 $L_2\geq |\nabla u|\geq L_1>0$ in $B(0, R)$. 

Using  Lemma  \ref{lemw}, we will prove that there exist $R>0$ and  $m>0$ such that for $K= \{ r\leq R, \sigma \in  S\}$ and $w=r^\gamma\psi(\sigma)$ then $u+ mw$  satisfies on $K$
\begin{equation}\label{234}
F(\nabla  (u+ mw), D^2 (u+m w) )-(L_f+1) (u+mw) > F(\nabla v, D^2 v) -(L_f+1) v 
\end{equation}  
with  $v \geq u+ mw$ on the boundary of $K$. 

A direct consequence of Lemma \ref{lemw} is that $\gamma>2$ for $a\neq A$. This implies that to conclude the proof we need the H\"older regularity of the derivative (see Proposition \ref{rglr1}) which gives, since $u$ is ${\cal C}^1$ and  
$f$ is H\"older's continuous :

    $$||u||_{{\cal C}^{2, \kappa}(\overline{\Omega})} \leq C_{A,f,\partial\Omega}. $$
for some  $ \kappa_{A,f}$.

 So we choose $\delta>0$ such that the sector $S=S_\delta $ is sufficiently close to the quarter of sphere, $\epsilon$ sufficiently small   and $a$ and $A$ sufficiently close to each other  in order that $\gamma_{a,A,\epsilon}^{S_ \delta }$ in  Lemma \ref{lemw} be such that $\gamma_{a,A,  \epsilon}^{S_\delta} \leq  2+ \kappa$, where 
 $\kappa$ is  
 recalled above.   In the following we drop for simplicity the index and exponents in $\gamma_{a,A}^S$ and  use the notation $\gamma$. 
 
 Referring to the notation of Lemma  \ref{lemw}, we now choose $R$  small enough in order that  $K\subset \Omega$, $C_{a,A, S_\delta , \gamma}  R^{\gamma-1} \leq \inf(  {L_1\over 2}, {L_2\over 2})$, 
and  such that- assuming  from now on and for simplicity that $\alpha <0$, the changes to bring for $\alpha >0$ being immediate -

$$\left({3L_2\over 2}\right)^\alpha R^{-2}\varepsilon >  2^{|\alpha|} L_1^\alpha C_{a,A, S_\delta , \gamma}  R^{-1} |f(u)|_\infty+ (L_f+1)$$
 Let  finally  $m<1$ be such that $(v-u)(R, \theta) \geq m R^\gamma \psi (\theta)$.   This is possible using the strict comparison principle in the following way:

  Let us  observe that  if $Q\in \Omega^{\star}\cap H_o\cap B_R(P)$ is such that $(v-u)(Q)=0$, $\partial_{x_1} (u-v) (Q)<0$. Indeed, let $B$ be a ball tangent to $\partial \Omega$ on $Q$, $B\subset \Omega\cap B_R(P)$.  By Theorem \ref{chopf},  $u < v$ in $B$ implies   $\partial_{x_1} (u-v) (Q)<0$.  Then by the continuity of $\partial x_1 (u-v)$ , there exists some neighborhood $V_Q $  and some $\delta_Q$ such that $(v-u)\geq \delta_Q x_1$ on $V_Q$. If $(v-u)(Q) >0$ the  same result is obvious by continuity.   Using a finite recovering  of  the sphere of center $P$ and radius $R$ by such neighborhoods     one gets 
    $$(v-u)(R, \sigma) \geq \delta x_1 \geq   \delta R d(x, \partial S) \geq {R \delta \over lip \psi} \psi (\sigma)\equiv m R^\gamma  \psi (\sigma).$$

We now observe that  by the choice of $m$ and $R$,  
${L_1\over 2}\leq |\nabla (u+ mw)|\leq  {3L_2\over 2}$ and then 
$\left({3L_3\over 2}\right)^\alpha \leq |\nabla (u+mw)|^\alpha \leq \left({L_1\over 2}\right)^\alpha$.  Also : 
\begin{eqnarray*}
|\nabla (u+ mw)|^\alpha {\cal M}_{a,A} (D^2u+ m  D^2 w) &\geq &|\nabla u|^\alpha {\cal M}_{a,A}(D^2 u)-m|\nabla w|_\infty L_1 2^{-\alpha} |f(u)|_\infty \\
&+& \left({3L_2\over 2}\right)^\alpha {\cal M}_{a,A}^- (m D^2 w)\\
&\geq& |\nabla u|^\alpha {\cal M}_{a,A}(D^2 u)+( L_f+1) w.
\end{eqnarray*}
Consequently one has 

 \begin{eqnarray*}
 |\nabla (u+ mw)|^\alpha {\cal M}_{a,A} (D^2u+ m  D^2 w) & -& (L_f+1) (u+m w) \\
& \geq&  |\nabla u|^\alpha{\cal M}_{a,A} (D^2 u)-  (L_f+1) u\\
  &\geq &   |\nabla v|^\alpha{\cal M}_{a,A} (D^2 v)- (L_f+1) v.
  \end{eqnarray*}
By Theorem \ref{comp} one derives that 
   $v \geq u+m w$. 
   
 \noindent    Suppose now that $\partial_{\vec \nu } u(P) = \partial_{\vec \nu} v(P)$, and   $\partial_{\vec \nu }^2 u(P) = \partial_{\vec \nu}^2 v(P)$. 
      This implies since $u$ and $v$ are in ${\cal C}^{2, \kappa}$ for the $\kappa$ given in Proposition \ref{rglr1}, that there exists some constant $c$ such that for all $r< R$, 
        $$(v-u)(r, \vec \nu) \leq c r^{2+\kappa}. $$
     This is a contradiction with 
    $v\geq u+ r^{\gamma} \psi (\vec \nu)$, since  $\vec \nu$  belongs to  $S_{\delta}$ as soon as $\delta$ is small enough.

\bigskip   
{\em Proof of Lemma \ref{lemw} }:  
In this proof, we need to compute a second order fully nonlinear operator on functions defined on the unit sphere. Since we shall use theories that have been developed only for fully nonlinear operators on functions on $\R^N$ we use an explicit system of coordinates on the sphere which  is  easy to manipulate and convenient for what we intend to  prove, but, of course, other choices are possible. 
 We denote by $\Sigma$ the  homeomorphism  which sends  $]0, {\pi\over 2}[\times ]-{\pi\over 2}, {\pi\over 2}[^{N-2}$  into  $S^+$,  defined by  
 
 $$ \sigma = \Sigma (\theta_1, \theta_2, \cdots, \theta_{N-1})$$
where $$\theta_1=\mbox{arctg}\frac{x_2}{x_1}\quad\mbox{and}\quad\theta_i=\mbox{arctg}\frac{x_{i+1}}{r_i},$$
here $r_i=\sqrt{\sum_{k=1}^{i} x_k^2}$ and, in the following, we shall use the following notations:
$$x'_i=(x_1,\dots,x_i) \quad\mbox{and}\quad(x'_i,0)=(x'_i,0,\dots,0).$$
It is easy to see that   if  $J = \left(\begin{array}{cc}
  0&-1\\
  1&0\end{array}\right)$:
$$\nabla \theta_1= (\frac{Jx'_2}{r_2^2},0)\quad
\mbox{and}
\quad \nabla \theta_i=(-\frac{x'_ix_{i+1}}{r_ir_{i+1}^2}, \frac{r_i}{r_{i+1}^2},0).$$
Hence 

$$\nabla \theta_i\cdot\nabla \theta_j=\delta_{ij}\frac{1}{r_{i+1}^2}.$$
Let $\psi$ be  a ${\cal C}^2$ function defined on  $]0, {\pi\over 2}[\times ]-{\pi\over 2}, {\pi\over 2}[^{N-2}$  and $w= r^\gamma \psi (\Sigma( \theta))$.  In the following, for simplicity,  we replace $\psi \circ \Sigma$ by $\psi$. 

 \noindent Obvious direct calculations give : 
$$\nabla  w= \gamma r^{\gamma-2} \psi x+ r^{\gamma}\psi_{\theta_i} \nabla \theta_i, $$
 \begin{eqnarray*}
D^2 w &= &r^{\gamma-2} \left( \psi_{\theta_i\theta_j}r^2(\nabla \theta_i\otimes\nabla \theta_j )\right. + \\
 &&+\gamma\psi_{\theta_i}(x\otimes\nabla  \theta_i+\nabla \theta_i\otimes x) \\
    && + \psi_{\theta_i}r^2(D^2\theta_i)+\\
    &&\left. +\gamma \psi\left(I+ (\gamma-2) {x\over r}\otimes {x\over r}\right)\right),
\end{eqnarray*}
summing over repeated indices.

We shall now compute the eigenvalues of each one of these matrices.
For $\Theta_i=x\otimes\nabla  \theta_i+\nabla \theta_i\otimes x$ it is easy to see that
$$ \Theta_i (x)=r^2 \nabla\theta_i\quad\mbox{and}\quad \Theta_i(\nabla  \theta_i)=|\nabla  \theta_i|^2x=\frac{1}{r_{i+1}^2}x.$$
This gives that the non trivial eigenvalues of $\Theta_i$ are 
$$\lambda_1(\Theta_i)=-\lambda_2(\Theta_i)=\frac{r}{r_{i+1}}.$$
Observe that $D^2\theta_1$ is zero outside a $2\times 2$ matrix and precisely for
$$J=\left(\begin{array}{ccc}
0 & -1 &0 \\
1 & 0 &0 \\
0& 0& 0
\end{array}\right)
 $$
 we get that
 $$
D^2\theta_1=\frac{1}{r_2^2}\left(J-2\frac{x^\prime_2}{r_2}\otimes\frac{Jx^\prime_2}{r_2}\right).
$$
whose non zero eigenvalues are $\lambda_1=\frac{1}{r_2^2}=-\lambda_2$.
Similarly, for $i>1$, $D^2\theta_i$ is a matrix which is zero outside of an $(i+1)\times ( i+1)$ matrix 
$$D_i=\left(\begin{array}{cc}
M_i & \frac{x^2_{i+1}-r^2}{r_ir_{i+1}^4} x^\prime_i \\
 \frac{(x^2_{i+1}-r^2)}{r_ir_{i+1}^4} x^\prime_i & -\frac{2r_ix_{i+1}}{r_{i+1}^4}
 \end{array}\right)
 $$
 with 
$$M_i=\frac{x_{i+1}}{r_ir_{i+1}^2}\left(-I+(x_i^\prime\otimes x_i^\prime)(\frac{1}{r_i^2}+\frac{2}{r_{i+1}^2})\right).
$$
Observe that $M_i$ (and $D_i$) has $i-1$ eigenvalues  for some  eigenvectors orthogonal to $x_i^\prime$:
$$\lambda(M_i)=-\frac{x_{i+1}}{r_ir_{i+1}^2}.$$
Since $M_i(x^\prime_i)=\frac{2x_{i+1}r_i}{r_{i+1}^4} x^\prime_i$
the other 2 eigenvalues of $D_i$ are given by the eigenvalues of the matrix
$$
\frac{1}{r_{i+1}^4}\left(\begin{array}{cc}
2x_{i+1}r_i & \frac{x_{i+1}^2-r_i^2}{r_i}\\
r_i(x_{i+1}^2-r_i^2) &-2r_ix_{i+1}
 \end{array}\right)
 $$
which are 
$$ \lambda_1=-\lambda_2=\frac{x_{i+1}^2+r_i^2}{r_{i+1}^4}=\frac{1}{r_{i+1}^2}.$$
Hence, for 
$$\Gamma(\sigma)=\left(\begin{array}{ccr}
\frac{r}{r_2} & 0 &\dots\\
0 & \frac{r}{r_3} & 0 \dots\\
\dots\\
0& \dots & 1
\end{array}\right),
 $$
 by writing the matrix  of $D^2 \psi$ in the basis generated by $(\nabla \theta_1,\cdots \nabla \theta_{N-1})$,  we denote with an abuse of notation, $\psi_{\theta_i\theta_j}r^2(\nabla \theta_i\otimes\nabla \theta_j )= \Gamma(\sigma)D^2(\psi) \Gamma (\sigma)$.

Using the properties of the operator ${\cal M}_{a,A}^-$ and using the same notation for the Pucci's operators on matrices $N\times N$ and  $(N-1)\times (N-1)$

\begin{eqnarray*}
{\cal M}_{a,A}^-(D^2 w)&\geq&   r^{\gamma-2} \left( {\cal M}_{a,A}^-(\Gamma(\sigma) D^2(\psi)\Gamma(\sigma))+ \right.\\
&&+(a-A)|\psi_{\theta_1}|(\frac{r^2}{r_2^2}+\gamma\frac{r}{r_{2}})\\
&& +(a-A)\sum_{i=1}^{N-1}|\psi_{\theta_i}| (\gamma\frac{r}{r_{i+1}}+ \frac{r^2}{r_{i+1}^2})+ \\
&&\left. +\sum_{i=1}^{N-1}\varepsilon_{a,A}(\psi_{\theta_i})\psi_{\theta_i}\frac{x_{i+1}r^2}{r_ir_{i+1}^2}+a\gamma (N+\gamma-2)\psi \right),
\end{eqnarray*}
where $\varepsilon_{a,A} (t) = a\ {\rm sign}^+ \ t + A\ {\rm sign}^- \ t$. 
Let  $H^\gamma_{a,A}$ the operator defined on  $S$ as 
\begin{eqnarray*}
H^\gamma_{a,A}(\psi)&:=&   {\cal M}_{a,A}^-(\Gamma(\sigma) D^2(\psi) \Gamma(\sigma))+ (a-A)|\psi_{\theta_1}|\ (\frac{r^2}{r_2^2}+\gamma\frac{r}{r_{2}})\\
&& +(a-A)\sum_{i=1}^{N-1}|\psi_{\theta_i}| (\gamma\frac{r}{r_{i+1}}+ \frac{r^2}{r_{i+1}^2})+ \\
&&+\sum_{i=1}^{N-1}\varepsilon_{a,A}(\psi_{\theta_i})\psi_{\theta_i}\frac{x_{i+1}r^2}{r_ir_{i+1}^2}.
\end{eqnarray*}
One can note that the coefficients in the definition of $H^\gamma_{a,A}$ are bounded on
every set $S_{\delta}$.    Furthermore  since $\Gamma (\sigma)$ is invertible,  
then the operator $H^\gamma_{a,A}$ is uniformly elliptic.    
Let us recall the definition  
 $$\overline{\lambda} (H_{a,A}^\gamma,S) = \sup\{ \lambda, \ \exists \ \psi >0 \ {\rm in} \ S ,\ H_{a,A}^\gamma (\psi) + \lambda \psi \leq \ 0 \ {\rm in} \ S\}.$$
 Using the results in \cite {IY},  $\overline{\lambda}$ is well defined and there exists $\overline{\psi} >0$ in $S$ such that 
 
 $$H_{a,A}^\gamma (\overline{\psi})+   \overline{\lambda} (H_{a,A}^\gamma,S)  \overline{\psi} = 0$$
  and $\overline{\psi} = 0$ on $\partial S$. Furthermore $\overline{\psi}$ is Lipschitz continuous.

\noindent    To conclude it is then sufficient to impose to $\gamma_{a,A, \epsilon}^S$ to be such that 
   $a\gamma_{a,A, \epsilon}^S (\gamma_{a,A, \epsilon}^S+ N-2) >  (\epsilon +  \overline{\lambda} (H_{a,A}^\gamma,S) )$ to get the desired result, and to take   $\psi$  as  some principal  positive eigenfunction  for the operator $H_{a,A}^\gamma$. 
   
   \bigskip
   
    We now prove the on $\gamma_{a,A, S_\delta}$:  
  
\noindent    Observe  that $H_{A,A}^\gamma (\psi) = A\Delta_S \psi$ where $\Delta_S$ 
   is the Laplace Beltrami operator on the sphere, and 
   $\overline{\lambda} (H_{A,A}^\gamma , S^+)= 2NA $. Then using classical elliptic estimates,  for all $\epsilon >0$,
    there exists $\delta>0$ such that for $|a-A|<\delta$,  there exists $\gamma\in (2, 2+\kappa_{A,f}),$
$$|\overline{\lambda} (H_{a,A}^\gamma,S_{\delta}) -2NA|\leq 2\epsilon. $$
This implies the result.

\bigskip

{\em Proof of Lemma \ref{lem1}.}
Suppose that $\phi$ is a ${\cal C}^2$ function, such that in a neighborhood of $P$, 
 $\partial \Omega$ coincides with  the graph $ x_N=\phi(x_1, \cdots x_{N-1})$. Without loss of generality we can suppose that $P=0$ and $e_N$ is normal to $\partial\Omega$ in $0$, hence  $\grad\phi(0)=0$.
 The Neumann boundary condition gives
 \begin{equation}\label{(2)}
 \partial_N u -\sum_{k=1}^{k= N-1}\partial_k u \partial_k\phi =- c (1+|\nabla \phi|^2)^{1\over 2};
 \end{equation}
 this together with the Dirichlet condition gives
 for $1\leq i\leq N-1$:
 \begin{equation}\label{(1)}
 (\partial_i u + \partial_N u \partial_i\phi)(x_1, \cdots, x_{N-1}, \phi(x_1,\cdots, x_{N-1}))=0.
 \end{equation}
Hence  $\partial_N u (0)  = -c$ and 
$\partial_i u(0)= 0$.

 Taking the derivative with respect to $x_j$ of (\ref{(1)}) with $j=1,\dots, N-1$ gives
 $$\partial_{ij} u(0)-c\partial_{ij} \phi(0)=0.$$  
Taking the derivative with respect to $j\in [1, N-1]$  of the identity (\ref{(2)})
 gives, 
 $$\partial_{Nj} u(0)=0.$$
 Finally
 $$D^2 u(0) = \left(\begin{array} {cc}
                    c D^2 \phi (0)&0\\
                     0& \partial_{NN}u (0)\end{array}\right)$$
 and then, by passing to the limit on the boundary  in  the equation  
                     $$|\nabla u|^\alpha {\cal M}_{a,A} (D^2 u)+f(u) = 0$$
one obtains   
$$u_{NN}(0)=\beta\left(-{\cal M} _{a,A} (cD^2 \phi)(0)-|c|^{-\alpha} f(0)\right),$$
  here ${\cal M}_{a,A}$ is understood as acting on $(N-1)\times(N-1)$ matrices and  $\beta=\frac{1}{a}$ or $\frac{1}{A}$ depending on the sign of $-{\cal M}_{a,A}  (cD^2 \phi)(0)-|c|^{-\alpha} f(0)$.

    \end{document}